\documentclass[12pt]{amsart}

\usepackage{enumerate}
 
\usepackage[active]{srcltx}
\usepackage{nicefrac}

\usepackage{amsthm,amssymb,setspace}
\usepackage[pagebackref,colorlinks,linkcolor=red,citecolor=blue,urlcolor=blue,hypertexnames=true]{hyperref}
\usepackage{amsrefs}
\usepackage[matrix, arrow]{xy}

\newcommand{\N}{\mathbb N}

\newcommand{\C}{\mathbb C}

\newcommand{\F}{\mathbb F}

\newcommand{\Z}{\mathbb Z}

\theoremstyle{plain}
\newtheorem{theorem}{Theorem}[section]
\newtheorem{corollary}[theorem]{Corollary}
\newtheorem{lemma}[theorem]{Lemma}
\newtheorem{proposition}[theorem]{Proposition}

\theoremstyle{definition}

\theoremstyle{remark}

\newtheorem{remark}[theorem]{Remark}

\begin{document}

\onehalfspace

\title{Convergent sequences in discrete groups}

\author{Andreas Thom}
\address{Andreas Thom, Universit\"at Leipzig, Germany}
\email{thom@math.uni-leipzig.de}

\begin{abstract} We prove that a finitely generated group contains a sequence of non-trivial elements which converge to the identity in every compact homomorphic image if and only if the group is not virtually abelian. As a consequence we show that a finitely generated group satisfies Chu duality if and only if it is virtually abelian.
\end{abstract}

\maketitle

\tableofcontents

\section{Introduction}

Let $G$ be a finitely generated discrete group. We call a sequence $(g_n)_{n \in \N}$ of elements of $G$ \emph{Bohr convergent to the neutral element}, if $\pi(g_n) \to 1$ for every homomorphism $\pi \colon G \to K$ to a compact group. We say that such a sequence is \emph{strongly Bohr convergent} if it is Bohr convergent and moreover $\pi(g_n)$ converges uniformly to the identity in unitary representations of a fixed finite dimension.

The aim of this note is to prove the following theorem.

\begin{theorem} \label{main}
Let $G$ be a finitely generated discrete group. Then, the following conditions are equivalent:
\begin{enumerate}
\item There exists a sequence of non-trivial elements which is strongly Bohr convergent to the neutral element.
\item The group is not virtually abelian.
\end{enumerate}
\end{theorem}

The assertion of the theorem is trivial if $G$ does not have enough homomorphisms into compact groups to separate elements of $G$, i.e.\ we are implicitly assuming that $G$ is maximally almost periodic. Being maximally almost periodic is equivalent to having a family of separating finite-dimensional unitary representations or the existence of an injection into a compact group. Since $G$ is assumed to be finitely generated, this condition is also equivalent to being residually finite by a well-known result of Mal'cev, see \cite{MR925989} for an elementary account. Note that there is always a \emph{net} which is Bohr convergent to the neutral element. There is a natural way of phrasing the result in terms of the Bohr compactification $bG$ of $G$. There is an extensive literature about the Bohr topology on various classes of topological groups and the existence or non-existence of Bohr convergent sequences has been an import question over the years (see \cite{MR2273361,MR1858143,her} and the references therein). Our result answers several open questions from \cite{MR1858143}. In particular, we provide an example of a non-trivial convergent sequence in the Bohr compactification of $\F_2$, a question which was asked by various authors. As a consequence of the main theorem, we can also show that the free group on two generators (and in fact any finitely generated group which is not virtually abelian) does not satisfy Chu duality (see Section \ref{chusec}). This problem was posed by Hsin Chu in his foundational work on Chu duality, see \cite{MR0195988,MR2273361}.

The following corollary is a consequence of our result which shows that images of word-maps can be surprisingly small. 

\begin{corollary}
Let $n \in \N$ and $\varepsilon>0$. There exists $w \in \F_2 \setminus \{e\}$ such that $$\|1_n - w(u,v)\| \leq \varepsilon, \quad \forall u,v \in U(n).$$
\end{corollary}

Here, $w \colon U(n) \times U(n) \to U(n)$ denotes the natural map, which is given by evaluating $w$ on the unitaries $u,v \in U(n)$. A quantitative version of this corollary is obtained in Remark \ref{remseq}.

The proof of the main theorem falls into a case study of the following three different cases:
1) $G$ has a free subgroup,
2) $G$ has no free subgroup and is not virtually abelian, and
3) $G$ is virtually abelian.
We will provide rather different arguments that deal with these various cases. The first case is dealt with in Section \ref{free} providing an explicit convergent sequence in the free group on two generators. The second case has been studied by S.\ Hern{\'a}ndez in \cite{her} using the Tits alternative and the structure theory of compact Lie groups. Hern{\'a}ndez showed the existence of a non-trivial sequence which is Bohr convergent to the neutral element.
We want to provide in Section \ref{nfree} an elementary argument for the existence of a \emph{strongly} Bohr convergent sequence. The proof is inspired by the ideas in \cite{her} but relies crucially on an extension of Jordan's theorem (see Proposition \ref{prop2}). 

It is well known \cite{glick,lep} that in the third case, there cannot be any Bohr convergent sequences.
Indeed, let $G$ be a virtually abelian group and assume that there exists a sequence $(g_k)_{k \in \N}$ such that $\pi(g_k) \to 1$ in every finite dimensional unitary representation. Since $G$ is virtually abelian, there exists an extension
$$1 \to A \to G \to F \to 1$$
with $A$ abelian and $F$ finite. Since $F$ has a faithful finite-dimensional unitary representation, we may assume that $g_k \in A$ for every $k \in \N$. Denote by $\widehat A$ the Pontrjagin dual of $A$, equipped with the normalized Haar measure $\mu$.

\begin{lemma} Let $A$ be a discrete abelian group.
Let $(g_k)_{k \in \N}$ be a sequence of non-trivial elements of $A$. Then, there exists $\chi \in \widehat A$ such that $\chi(g_k) \not \to 1$.
\end{lemma}
\begin{proof}
Assume that $\chi(g_k) \to 1$ for all $\chi \in \widehat A$. We conclude from the Dominated Convergence Theorem that
$$0 = \lim_{k \to \infty} \int_{\widehat A} \chi(g_k) \ d\mu(\chi) = \int_{\widehat A} 1 \ d\mu(\chi) = 1.$$
This proves the claim.
\end{proof}

The character $\chi \colon A \to S^1$ provides a one-dimensional unitary representation. Since $A$ has finite index in $G$, we may induce $\chi$ to a finite-dimensional unitary representation $\pi={\rm Ind}_A^G(\chi)$ of $G$. Since $\chi \subset {\rm Res_A^G}\pi$, the unitary representation $\pi$ satisfies $\pi(g_k) \not\to 1$ as $k \to \infty$.

\section{Groups without free subgroups} \label{nfree}

In this section we deal with case 2), i.e.\ groups which are not virtually abelian and do not contain free subgroups. The proof is inspired by ideas appearing in \cite{her}.
Let us first recall basic facts about the unitary groups. Let $n \in \N$ be an integer and $U(n)$ be the group of unitary matrices. We endow $U(n)$ with the metric coming from the operator norm. More specifically, we set
$$\ell(g) := \|1-g\|, \quad \mbox{and} \quad d(g,h) := \ell(gh^{-1}) = \|g-h\|.$$
We use the notation
$[g,h] = ghg^{-1}h^{-1}$ and $\bar g = g^{-1}$.

\begin{lemma} \label{comp}
Let $g,h \in U(n)$. The following relations hold:
$$\ell(gh) \leq \ell(g) + \ell(h), \quad \ell(gh\bar g)=\ell(h) \quad \mbox{and} \quad \ell([g,h]) \leq 2 \cdot \ell(g) \ell(h).$$ 
\end{lemma}
\begin{proof}
The last inequality is the only non-trivial assertion. Let us compute:
$$\ell([g,h])=\|1 - gh\bar g \bar h\| = \|hg -gh\| = \|(1-h)(1-g) - (1-g)(1-h)\| \leq 2 \cdot \ell(g) \ell(h).$$
This finishes the proof.
\end{proof}

In the proof, we are using the following structure result about subgroups of unitary groups.
\begin{proposition} \label{prop1}
Let $G$ be a group without free subgroups, $n \in \N$, and let $\pi \colon G \to U(n)$ be a finite-dimensional unitary representation. Then, $\pi(G)$ is virtually abelian.
\end{proposition}
\begin{proof}
Since $\pi(G)$ is linear and does not contain a free subgroup, it follows from the Tits alternative \cite{MR0286898}, that $\pi(G)$ is virtually solvable. Let $H \subset \pi(G)$ be a solvable subgroup of finite index. The closure $\overline{H} \subset U(n)$ is a solvable compact Lie group. By Theorem $29.44$ in \cite{MR551496}, the connected component of $1$ in $\overline{H}$ is abelian. Since $\overline{H} \subset U(n)$ is a Lie group, it has only finitely many connected components and it follows that $\overline{H}$ is virtually abelian. This implies that $H$ and hence $\pi(G)$ are virtually abelian.
\end{proof}

A classical theorem of Camille Jordan \cite{jordan} says that for fixed $n \in \N$ there exists an integer $m(n) \in \N$ such that every finite subgroup of $GL(n,\C)$ has an abelian subgroup of index at most $m(n)$. Note that any finite subgroup is conjugate to a subgroup of $U(n)$. We will need the following generalization of Jordan's theorem.

\begin{proposition}\label{prop2}
Let $n \in \N$. There exists $m(n) \in \N$ such that every $G \subset U(n)$ without free subgroups has an abelian normal subgroup of index at most $m(n)$.
\end{proposition}
\begin{proof}
Proposition \ref{prop1} says that $G$ has to be virtually abelian.
We prove the statement by induction over $n$. If $n=1$, there is nothing to prove since $U(1)$ is abelian. 

Let $K \subset G$ be the subgroup with is generated by elements $g \in G$ for which $\ell(g) < n^{-1/2}/2.$ The group $K$ is a normal subgroup and there exists a universal constant $c(n)$ such that $[G:K] \leq c(n)$.
Let $A \subset K$ be a maximal abelian normal subgroup. 
If $A$ does not consist of multiples of identity, then the normalizer of $A$ is a finite extension of $U(n_1) \times \cdots \times U(n_k)$ for some $n_1,\dots,n_k$ such that $n_1 + \cdots + n_k \leq n$ and $k \geq 2$. This is seen by noting that $A$ can be simlutaneously diagonalized. Then, the normalizer is generated by the centralizer (i.e.\ a product of unitary groups) and those permutations, which permute blocks  of equal size consisting of identical eigenvalues.
In fact, the normalizer is a product of groups of the form $U(k)^{\times l} \rtimes S_l$ and the index of $U(n_1) \times \cdots \times U(n_k)$ in the normalizer is bounded by $e(n):=n!$. 

Hence, upon passing to a subgroup of index at most $e(n)$ in $K$, we may assume that, $K \subset U(n_1) \times \cdots \times U(n_k)$ and we can finish the proof using the induction hypothesis.
Indeed, let $\pi_i(K)$ be the image of $K$ under the projection to $U(n_i)$. By induction, there exists a normal abelian subgroup $A_i \subset \pi_i(K)$ of index at most $m(n_i)$. Since 
$K \subset \prod_{i=1}^n \pi_i(K)$
we conclude that $K$ admits an abelian subgroup of index at most $m({n_1}) \cdots m(n_k)$. Let $d(n)$ be the maximum of all values of $m(n_1) \cdots m(n_k)$ over all partitions of natural numbers less or equal to $n$. We conclude that $G$ has a normal abelian subgroup of index $m(n) := c(n) d(n)e(n)$.

Hence, we may assume in the continuation of the proof that $A$ acts with multiples of identity. If $A = K$ we are done so that we may also assume that $A \subsetneq K$. Since $A \subset K$ is a normal subgroup of finite index, there exist an element $g \in K$ which is not a multiple of identity such that $\min_{\lambda \in S^1} d(\lambda, g)$ is smallest possible. Since $K$ is not abelian, not all generators of $K$ can be multiples of identity and we have $\min_{\lambda \in S^1} d(\lambda, g) < n^{-1/2}/2$. Let $h$ be a generator of $K$ and compute 
\begin{eqnarray*}\ell([g,h]) &=& \|gh - hg\| \\
&=& \min_{\lambda \in S^1}\|(\lambda - g)(1-h) - (1-h)(\lambda - g)\| \\
&\leq& n^{-1/2} \cdot \min_{\lambda \in S^1} d(\lambda,g)\\
&<& n^{-1}/2.
\end{eqnarray*}
In particular, $\min_{\lambda \in S^1} d(\lambda, [g,h]) < \min_{\lambda \in S^1} d(\lambda, g)$ and hence $[g,h]$ is a multiple of identity. The determinant of $[g,h]$ equals one and hence $[g,h] = \lambda$ for some $\lambda \in S^1$ with $\lambda^n = 1$. Now, for an $n$-th root of unity either $\ell([g,h])=|\lambda -1| \geq n^{-1}/2$ or $\lambda =1$. We conclude that $[g,h]=e$ and thus $g$ commutes with all generators of $K$. We finally get that $g \in A$ by maximality of $A$. However, this contradicts our choice of $g \in K$ and finishes the proof. 
\end{proof}

The following result is of independent interest and will not be used in the sequel.
\begin{corollary}
Let $n \geq 2$. The set
$$\{(u,v) \in U(n) \times U(n) \mid \langle u,v \rangle \mbox{ does not contain a free subgroup}\}$$ is not dense in $U(n) \times U(n)$.
\end{corollary}
\begin{proof}
We first claim that there exists an integer $m$ such that
$[u^m,v^m] = e$ if $\langle u,v \rangle$ has no free subgroups.
We set $G := \langle u,v \rangle$. By the preceding theorem, there exists a normal abelian subgroup $A \subset G$ of index at most $m(n)$. Hence, $G/A$ has order at most $m(n)$. Putting $m:= m(n) !$, we obtain that $u^m$ and $v^m$ are trivial in $G/A$ and hence, $u^m,v^m \in A$. We conclude $[u^m,v^m] = e$.

By continuity of the multiplication, $[u^m,v^m] = e$ holds for all pairs $(u,v)$ in the closure of the set $\{(u,v) \in U(n) \times U(n) \mid \langle u,v \rangle \mbox{ does not contain a free subgroup}\}$. However, it is well-known that a generic pair of unitaries in $U(2)$ generates a free group. Indeed, the existence of free subgroups of $SO(3)$ (and hence its double cover $SU(2)$) was already known to F.\ Hausdorff \cite{MR1511802}. It follows that for each of the countably many non-trivial words, the set of pairs in $U(2) \times U(2)$ which satisfies the word is an {\it proper} real-algebraic subvariety. In particular, its Haar measure is zero. We conclude that the measure of the set of pairs which satisfy a non-trivial word is zero.
\end{proof}

\begin{proposition} \label{prop3}
Let $n \in \N$ and $G$ be a finitely generated group without free subgroups which is not virtually abelian. There exists $g \in G \setminus \{e\}$ such that $\phi(g) = 1$ for every homomorphism
$$\phi \colon G \to U(n).$$
\end{proposition}
\begin{proof}
From Proposition \ref{prop2} we conclude that every homomorphic image of $G$ in $U(n)$ admits a finite index subgroup of index at most $m(n)$ which is abelian. Let $K \subset G$ be the intersection of all normal subgroups of $G$ which are of index at most $m(n)$. Since $G$ is finitely generated, $K$ is again of finite index. Since $G$ is not virtually abelian, $K$ is not abelian and we may set
$g:= [u,v]$ for some pair $u,v$ of non-commuting elements in $K$. This finishes the proof. 
\end{proof}

\begin{corollary} \label{nonfree}
Let $G$ be a finitely generated  group without free subgroups which is not virtually abelian. Then, there exists a sequence $(g_k)_{k \in \N}$ of non-trivial elements such that
$\pi(g_k) \to 1$ for every homomorphism $\pi \colon G \to K$ into a compact group.
Moreover, if $\pi \colon G \to U(n)$ for some $n \in \N$, then $\pi(g_k) = 1$ for $k$ large enough.
\end{corollary}
\begin{proof}
By Proposition \ref{prop3}, there exists $g_k \in G$ such that $\pi(g_k) = 1$ for every homomorphism $\pi \colon G \to U(k)$.
Let $\pi \colon G \to K$ be a homomorphism into a compact group. We claim that $\pi(g_k) \to 1$ as $k \to \infty$. By the analysis leading to the Peter-Weyl theorem \cite{MR551496}, every  compact group has a separating family of finite-dimensional unitary representations $\sigma_{\alpha} \colon K \to U(n_{\alpha})$. In particular, there exists a continuous injection
$$\sigma \colon K \to \prod_{\alpha} U(n_{\alpha}).$$
Since $K$ is compact, $\pi$ is a homeomorphism onto its image. This reduces the argument to the case $K= U(n)$ for some fixed $n \in \N$. However, in this case $\pi(g_k) = 1$ for $k \geq n$. This proves the claim. 
\end{proof}

\section{Groups with free subgroups} \label{free}

In this section, we deal with the case that $G$ has a free subgroup. We will provide an explicit sequence of elements in the free group on two generators which is strongly Bohr convergent to the neutral element.

\begin{proposition} \label{firstmain}
There exists a sequence $(z_k)_{k \in \N}$ of non-trivial elements of the free group $\F_2$, such that for every $r \in \N$ there exists a constant $C(n)$ such that for every unitary representation $\pi \colon \F_2 \to U(n)$, $$\ell(\pi(z_k)) \leq C(n) \cdot 2^{-k}, \quad \forall k \geq 1.$$
\end{proposition}

Before we prove the proposition, we need to prove a lemma.

\begin{lemma}\label{clear}
Let $n \in \N$ be an integer and $\varepsilon>0$. There exists $q(n,\varepsilon) \in \N$ such that for all $g \in U(n)$, there exists
$m \in \{1,2,\dots,q(n,\varepsilon)\}$
with $\ell(g^m) \leq \varepsilon$.
\end{lemma}
\begin{proof}
Suppose that there exists a natural number $q \in \N$ such that $\ell(g^m) > \varepsilon$ for all $m \in \{1,2,\dots,q \}$. We conclude that $d(g^m,g^{m'}) = \ell(g^{|m-m'|})> \varepsilon$ for all $1 \leq m,m' \leq q$. Hence, the balls of radius $\varepsilon/2$ around $g^m$ are pairwise disjoint for $1 \leq m \leq q$. Denote the Haar measure on $U(n)$ by $\mu$. We conclude that
$\mu(U(n)) \geq q \cdot \mu( \{u \in U(n) \mid \ell(u) \leq \varepsilon/2 \})$
and hence 
$$q \leq \frac{\mu(U(n))}{\mu( \{u \in U(n) \mid \ell(u) \leq \varepsilon/2 \})}.$$
Hence, we may set $q(n,\varepsilon):= q$ for some $q > \frac{\mu(U(n))}{\mu( \{u \in U(n) \mid \ell(u) \leq \varepsilon/2 \})}$. This finishes the proof.
\end{proof}

\begin{proof}[Proof of Proposition \ref{firstmain}] Let $\F_2$ be generated by the letters $a$ and $b$. It is well-known, that the set $\{a^nba^{-n} \mid n \in \Z\} \subset \F_2$ is free, i.e.\ there are no non-trivial relations among the elements $a^nba^{-n}$. Assume now that $w \in \F_2$ together with the set $\{ w_{k,n} \mid (k,n) \in \Z \times \N\} \subset \F_2$ are free by choosing a bijection $\phi \colon \{\star\} \cup \Z \times \N \stackrel{\sim}{\to} \Z$ and setting $w_{k,n} = a^l b a^{-l}$ for $l= \phi(k,n)$ and $w = a^l b a^{-l}$ for $l = \phi(\star)$. 
We set
$$v_{k,1} = \begin{cases} w_{k,1} w^k \bar w_{k,1}, & \forall k \geq 1 \\ w_{k,1} w^{k-1} \bar w_{k,1}, & \forall k \leq 0\end{cases} .$$
and define by induction:
$$v_{k,n+1} := w_{k,n+1} \left[[v_{k-1,n},v_{k,n}], v_{k+1,n}\right] \bar w_{k,n+1}.$$
It is easy to see by induction that $v_{k,n} \neq e$ for all $(k,n) \in \Z \times \N$. Indeed, we claim that the set $\{v_{k,n} \mid k\in \Z\}$ is free for each $n \in \N$. The claim follows from the obvervation that $[[a,b],c] \neq e$ for a basis of a free group of rank $3$ and the fact that conjugating with $w_{k,n}$ in the definition of the $n$th stage of the double-sequence produces again a sequence forming a free subset.
 
Let $\phi \colon \F_2 \to U(r)$ be a homomorphism and consider the induced length function on $\F_2$, i.e.
$\ell(g) = \|1 - \phi(g)\|$ for all $g \in \F_2.$ Note that $\ell(g) \leq 2$ for all $g \in \F_2$.
We compute from Lemma \ref{comp}
\begin{equation} \label{eq1} \ell(v_{k,n+1}) \leq 2^2 \cdot \ell(v_{k-1,n}) \cdot \ell(v_{k,n}) \cdot \ell(v_{k+1,n}).\end{equation} Hence, using Equation \ref{eq1} at the spots $k-1,k$ and $k$, we get:
\begin{eqnarray} \label{eq2}
\ell(v_{k,n+2}) &\leq& 2^8 \cdot \ell(v_{k-2,n}) \cdot \ell(v_{k-1,n})^2 \cdot \ell(v_{k,n})^3 \cdot \ell(v_{k+1,n})^2 \cdot \ell(v_{k+2,n})
\end{eqnarray}
Using $\ell(g) \leq 2$, we get $\ell(v_{k,n+2}) \leq 2^{14} \cdot \ell(v_{k,n})^3.$
In particular, if $\ell(v_{k,n}) \leq 2^{-15}$ for some $k \in \Z$, then 
\begin{equation} \label{cruc} \ell(v_{k,n+2}) \leq 2^{-16} \cdot \ell(v_{k,n}).\end{equation}
Moreover, if $\ell(v_{k,n})\leq 2^{-15}$, then 
Equation \ref{eq2} and the trivial bound $\ell(g) \leq 2$ give:
$$\ell(v_{k\pm 1,n+2}) \leq 2^{15}  \cdot \ell(v_{k,n})^2 \leq 2^{-15}.$$

We see that the property $\ell(v_{k,n}) \leq 2^{-15}$ spreads with speed $1$ in the $k$-direction as we increase the index $n$ by $2$. Moreover, by Equation \ref{cruc}, once $\ell(v_{k,n}) \leq 2^{-15}$ holds, the quantity $\ell(v_{k,n})$ decays exponentially in $n$. 
By Lemma \ref{clear}, there exists $q:=q(r,2^{-15})$ such that there exists $m \in \{1,2,\dots,q \}$ with $\ell(w^m)\leq 2^{-15}$.

We set $z_n := v_{0,2n+1}$ for $n \in \N$ and claim that the sequence $(z_n)_{n \in \N}$ solves the problem. 
Indeed, we can conclude from our observations above that
$$\ell(v_{k,2n}) \leq 2^{-15}, \quad \forall (k,n) \in \N \times \Z \mbox{ with } q - n \leq k \leq n.$$
Hence $\ell(z_q) \leq 2^{-15}$ and Equation \ref{cruc} implies that
$\ell(z_n)$ converges to zero and we have the estimate
$$\ell(z_{n}) \leq 2^{-15 - 16(n- q(r,2^{-15}))}, \quad \forall k \in \N \colon k \geq  2q.$$ 
This proves the claim.
\end{proof}

One can also put the statement of Proposition \ref{firstmain} in terms of so-called word-maps. Note that every $w \in \F_2$ gives rise to a continuous evaluation map $w \colon U(n) \times U(n) \to U(n)$ which is called the \emph{word-map} associated with $w$. We can state the following two immediate corollaries.

\begin{corollary} \label{seq}
Let $n \in \N$ and $\varepsilon>0$. There exists a non-trivial element $w \in \F_2$, such that the natural map
$w \colon U(n) \times U(n) \to U(n)$
satisfies $$\ell(w(g,h)) \leq \varepsilon,$$ for all $g,h \in U(n)$.
\end{corollary}

\begin{corollary} \label{cor2}
Let $n \in \N$ and $\varepsilon>0$. There exist $k \in \N$, unitaries $u_1,\dots,u_k$ and signs $\varepsilon_1,\dots,\varepsilon_k \in \{\pm 1\}$, such that the continuous map $\phi \colon U(n) \to U(n)$ which is defined as
$$\phi(v) := u_1 v^{\varepsilon_1} \cdots u_k v^{\varepsilon_k} \in U(n)$$
is non-trivial and satisfies $\ell(\phi(v)) \leq \varepsilon$ for all $v \in U(n)$.
\end{corollary}

\begin{remark}
Murray Gerstenhaber and Oscar Rothaus showed in their work on solvability of equations over groups \cite{MR0166296} that continuous maps 
$\phi(v) = u_1 v^{\varepsilon_1} \cdots u_k v^{\varepsilon_k}$
as in Corollary \ref{cor2} are surjective if $\sum_{i=1}^k \varepsilon_i \neq 0$.
The preceding corollary shows how drastically this conclusion can fail without any assumption on the sum of the exponents.
\end{remark}

\begin{remark} \label{remseq}
Consider the free group $\F_2$ with generators $a,b$. We denote the word-length of a word $w \in \F_2$ with respect to the set $\{a,\bar a, b, \bar b\}$ by $L(w) \in \N$. 

There is a quantitative version of Corollary \ref{seq} which says that there exists a sequence $(w_k)_{k \in \N}$ of non-trivial elements in $\F_2$, constants $\alpha>0$ and $c_1(n)>0$, such that $L(w_k) \to \infty$ as $k \to \infty$ and
$$ \ell(w_k(u,v)) \leq \exp\left( - c_1(n) \cdot L(w_k)^{\alpha}\right), \quad \forall u,v \in U(n).$$
The construction above yields $\alpha = \log_{10} 3 - \varepsilon$ for every $\varepsilon>0$. A more refined procedure yields $\alpha = \log_{14} 4 - \varepsilon,$
for every $\varepsilon>0$. Note that $\log_{14} 4 = 0,525\overline{2}>1/2.$

This has to be contrasted by a result of V.\ Kaloshin and I. Rodnianski in \cite{MR1873135}. They prove that for almost all pairs $u,v \in SO(3)$, there exists a constant $c_2(u,v)>0$ such that one has
$$\ell(w(u,v)) \geq \exp(- c_2(u,v) \cdot L(w)^2 ).$$ A.\ Gamburd, D.\ Jakobson and P.\ Sarnak conjectured in \cite{MR1677685} that a similar inequality holds for almost all $u,v \in SO(3)$ with an exponent which is linear in $L(w)$.
\end{remark}

Let us now come back to the proof of Theorem \ref{main}.

\begin{proposition} \label{freepro}
There exists a sequence $(z_k)_{k \in \N}$ of non-trivial elements of the free group $\F_2$, such that for every homomorphism $\pi \colon \F_2  \to K$ to a compact group $$\pi(z_k) \to 1, \quad \mbox{as} \quad k \to \infty.$$
Moreover, the convergence is uniform when restricted to unitary representations of a fixed dimension.
\end{proposition}
\begin{proof} First of all, we may perform
a reduction to finite-dimensional unitary representations as in the proof of Proposition \ref{nonfree}.
Now, Proposition \ref{firstmain} implies the claim. 
\end{proof}

The proof of Theorem \ref{main} follows by combining Corollary \ref{nonfree}, Proposition \ref{freepro} and the remarks about virtually abelian groups at the end of the Introduction.

\section{Chu duality} \label{chusec}

In \cite{MR0195988}, Hsin Chu studied a concept of duality for maximally almost periodic groups which ought to generalize both Pontrjagin duality for (locally compact) abelian group and Tannaka-Krein duality for compact groups. Let $G$ be a locally compact group and set
$$G^{x}_n := \{ \phi \colon G \to U(n) \mid \phi \mbox{ a continuous homomorphism} \}.$$
We endow $G^{x}_n$ with the compact-open topology. Note that if $G$ is discrete, $G^x_n$ is compact. The set $G^x := \coprod_{n \geq 1} G^x_n$ comes equipped with natural continuous operations
$$\oplus \colon G^x_n \times G^x_m \to G^x_{n+m} \quad \mbox{and} \quad \otimes \colon G^x_n \times G^x_m \to G^x_{nm}$$
which are inherited from the corresponding operations on $U := \coprod_{n \geq 1} U(n)$. Moreover, for $\phi \in G^x_n$ and $u \in U(n)$, we may consider $u \phi u^* \in G_x^n$.

Chu proceeds by defining
$G^{xx}$ to be the set of continuous degree-preserving maps $\mu \colon G^x \to U$ such that
$\mu( \phi \oplus \phi') = \mu(\phi) \oplus \mu(\phi')$, $\mu(\phi \otimes \phi') = \mu( \phi) \otimes \mu(\phi')$ for all $\phi,\phi' \in G^x$ and $\mu(u \phi u^*) = u \mu(\phi) u^*$ for all $\phi \in G^x$ and $u \in U$ of the same degree. Moreover, he endows the set $G^{xx}$ with the compact-open topology and shows that the multiplication in $U$ induces a natural group structure on $G^{xx}$ which is compatible with this topology. The topological group $G^{xx}$ is nowadays called the \emph{Chu dual} of $G$. The evaluation map
${\rm ev} \colon G \to G^{xx}$ which is given by $g \mapsto \{ \phi \mapsto \phi(g) \}$ defines a natural continuous homomorphism. A group is said to satisfy \emph{Chu duality} if and only if ${\rm ev} \colon G \to G^{xx}$ is a homeomorphism. Note that the Pontrjagin duality theorem says that locally compact abelian groups satisfy Chu duality, the Tannaka-Krein duality theorem is the corresponding assertion for compact groups. Chu ended his paper \cite{MR0195988} with the question whether the free group on two generators satisfies his notion of duality.

As a corollary of Theorem \ref{main}, we can show:
\begin{corollary} \label{chu}
If a finitely generated group satisfies Chu duality, then it is virtually abelian. In particular, the free group on two generators does not satisfy Chu duality.
\end{corollary}
\begin{proof} Let us prove that any finitely generated group which is not virtually abelian does not satisfy Chu duality by explicitly showing that ${\rm ev} \colon G \to G^{xx}$ cannot be a homeomorphism.

Since $G$ is finitely generated and not virtually abelian, there exists a non-trivial sequence $(g_k)_{k \in \N}$ which is strongly Bohr convergent to the neutral element in $G$. This implies that ${\rm ev}(g_n)$ converges to $1 \in G^{xx}$ in the compact-open topology. Indeed, any compact subset of $G^x$ is contained in a finite union $\coprod_{1 \leq n \leq l}G^x_n$ and convergence was uniform on unitary representations of a fixed dimension. We conclude that the topology on $G^{xx}$ is not discrete. Hence, ${\rm ev} \colon G \to G^{xx}$ cannot be a homeomorphism.
\end{proof}


\begin{remark}
The converse to Corollary \ref{chu} is well-known. Indeed, if $G$ is virtually abelian, then there is a universal bound on the dimension of an irreducible representation and this simplifies the study of $G^{xx}$ considerably.
\end{remark}

\begin{remark}
If one bases the unitary duality theory on representations on Hilbert spaces of finite \emph{and} infinite dimensions, one obtains a well-working duality theory, i.e.\ the canonical homomorphism into the analogously defined bi-dual is an isomorphism of topological groups if the group is locally compact.
\end{remark}

\begin{remark} Let $G$ be a finitely generated group which is not virtually abelian.
Since $G$ is countable and discrete, $G^x_n$ is separable and compact for all $n \in \N$. This implies that $G^{xx}$ carries the structure of a second countable complete metric space. We conclude that the homomorphism ${\rm ev} \colon G \to G^{xx}$ cannot be onto, since that would imply that ${\rm ev} \colon G \to G^{xx}$ is a homeomorphism; see \cite{MR0364542} for details. Hence, $G$ is not  Chu semi-reflexive. We also note that ${\rm ev}(G)$ cannot be complete.
\end{remark}

\section*{Acknowledgment}

I want to thank Beno\^{i}t Collins for mentioning the problem during a visit at the University of Ottawa in 2007 and Marc Burger for a stimulating discussion about Jordan's theorem.

\end{document}